\documentstyle[12pt,leqno,amssymb,graphics]{amsart}

\newtheorem{thm}{Theorem}[section]
\newtheorem{lemma}[thm]{Lemma}
\newtheorem{prop}[thm]{Proposition}

\theoremstyle{definition}
\newtheorem{dfn}[thm]{Definition}
\theoremstyle{remark}
\newtheorem{remark}[thm]{Remark}
\newtheorem{ques}[thm]{Question}

\begin{document}

\newcommand{\pr}{\protect\ref}
\newcommand{\su}{\subseteq}
\newcommand{\pa}{{\partial}}
\newcommand{\C}{{\mathcal C}}
\newcommand{\A}{{\mathcal A}}
\newcommand{\R}{{\Bbb R}}
\newcommand{\B}{{\Bbb B}}
\newcommand{\E}{{{\Bbb R}^3}}
\newcommand{\G}{{\Bbb G}}
\newcommand{\Z}{{\Bbb Z}}
\newcommand{\im}{{Imm(F,\E)}}
\newcommand{\eq}{{\equiv}}
\newcommand{\e}{{\epsilon}}
\newcommand{\p}{{p_1,\dots,p_n}}
\newcommand{\lm}{{\lambda}}
\newcommand{\dd}{{\dots}}

\newcounter{numb}

\title{Order One Invariants of Immersions}
\author{Tahl Nowik}
\address{Department of Mathematics and Computer Science, Bar-Ilan University, Ramat-Gan 52900,  Israel.}
\email{tahl@@macs.biu.ac.il}
\date{March 25, 2001}

\begin{abstract}
We classify all order one invariants of immersions
of a closed orientable surface $F$ into $\E$, with values in an arbitrary
Abelian group $\G$.
We show that for any $F$ and $\G$ and any regular homotopy class $\A$ of immersions of $F$ into $\E$,
the group of all order one invariants on $\A$ is isomorphic to
$\G^{\aleph_0} \oplus \B \oplus \B$ where $\G^{\aleph_0}$ is the group of all functions from a set of cardinality $\aleph_0$ into $\G$ and
$\B=\{ x\in\G  :  2x=0 \}$.

Our work includes foundations for the study of finite order 
invariants of immersions of a closed orientable surface into $\E$, 
analogous to chord diagrams and the
1-term and 4-term relations of knot theory.
\end{abstract}

\maketitle

\section{Introduction}\label{0}

The notion of finite order invariants has developed in knot theory.
We extend it here to the setting of immersions of a closed orientable surface 
into $\E$. We give a general analysis of finite order invariants, finally classifying all order one invariants (Theorem \pr{t1}).
A small subclass of order one invariants has previously been studied, namely, the ``local invariants'', studied in [G], [N1], [N3].

The structure of the paper is as follows:
In Section \pr{A}  we describe the self-intersection appearing in immersions which correspond to the codimension 1 (abbreviated codim 1) strata in the space of immersions. They will be the basis of our definition of finite order invariants. 
In Section \pr{B} we show that there is no continuous choice of co-orientation
for the codim 1 strata in the space of immersions. We show that this is not an obstacle for defining finite order invariants; we construct the combinatorial object on which order $n$ invariants are defined up to order $n-1$ invariants.
In Section \pr{C} we analyze the relations on invariants arising from local 
2-parameter families of immersions. 
In Section \pr{D} we classify all order one invariants.
In Section \pr{E} we show that for $\G={\Bbb Z}/2$, 
the analogous classification does not hold for invariants of order $n>1$;
in that we also complete the analysis of [N2].

\section{Self intersection}\label{A}

The space of all immersions of a closed orientable surface $F$ into $\E$ is naturally stratified by the types of self intersection appearing in the immersions; see [HK] for 
local analysis. Parts of our work will also be local;
given an immersion $i:F\to\E$ and a self intersection point of $i$ we will look at all deformations of $i$ which move $F$ only in a small neighborhood of that intersection point. The \emph{local stratification} will then mean the stratification of this smaller space of immersion and \emph{local codimension} of the local strata will again refer to this smaller space.
We will take a close look at the local stratification up to codim 2.
We discuss the codim 0 and 1 strata here. 
The codim 2 strata will be discussed in Section \pr{C}. 

The codim 0 strata corresponds to stable self intersection, which are double lines and triple points. 
The codim 1 strata divide into four types 
which (following [G]) we call: $E$, $H$, $T$, $Q$. 
In the notation of [HK] they are respectively
$A_0^2|A_1^+$, $A_0^2|A_1^-$, $A_0^3|A_1$, $A_0^4$.
The four types may be demonstrated by the following local representatives, where formulae
in 3-space defining the different sheets involved in the self intersection, are given. A representative of the  codim 1 strata is obtained from the formulae below by setting $\lm=0$. Letting $\lm$ vary, we obtain a 1 parameter family of immersions which is transverse to this codim 1 stratum.  

$E$: \ \ $z=0$, \ \ $z=x^2+y^2+\lm$. See Figure \pr{fet}, ignoring the vertical plane.

$H$: \ \ $z=0$, \ \ $z=x^2-y^2+\lm$. See Figure \pr{fht}, ignoring the vertical plane.

$T$: \ \ $z=0$, \ \ $y=0$, \ \ $z=y+x^2+\lm$. See Figure \pr{ftq}, ignoring the vertical plane.

$Q$: \ \ $z=0$, \ \ $y=0$, \ \ $x=0$, \ \ $z=x+y+\lm$. This is simply four planes passing through one point, no three of which intersect in a line

A self intersection of local codim 1 (namely $E$, $H$, $T$ or $Q$) will be called a CE point (after the first and last letters of  ``co-dimension one'').  By CE point we will also refer to the point in $\E$  where this self-intersection takes place.
Once an orientation is chosen for the surface, the above four types of self intersection split into twelve types, as we explain below.
A choice of one of the two sides of the local codim 1 stratum at a given point of the stratum, is represented by the choice of $\lm<0$ or $\lm>0$ in the formulae above. We will refer to such a choice as a 
\emph{co-orientation for the configuration of the self intersection.} 
A completely different notion of co-orientation that we will encounter is the following: 
If $i:F\to\E$ is an immersion and $F$ is oriented then its orientation and that of $\E$ induce a co-orientation for $i(F)$ in $\E$. 
To avoid confusion between the two notions,
we will use the term ``co-orientation'' only for the former.
For the latter we will speak of the ``preferred side of $i(F)$ in $\E$''. 

We now present the twelve types of self intersection in the oriented setting and specify a co-orientation wherever possible:

Type $E$: the configuration of the self intersection at the two sides of the stratum is distinct even with no orientation on the surface, namely, for $\lm<0$ there is an additional 2-sphere in the image of the immersion, and we choose this side ($\lm<0$) as our positive side for the co-orientation. Now, this additional 2-sphere is made of two 2-cells; if the surface is oriented then we may distinguish three different types of $E$ self intersections, which we denote by $E^a$ where $0\leq a \leq 2$ denotes the number of such 2-cells
for which the little 3-cell bounded by the 2-sphere lies on their non-preferred side in $\E$. 

Type $H$: the configuration of the self intersection on the two sides of this stratum are indistinguishable with no orientations. If the surface is oriented
then we may compare the orientations of the two sheets at the time of tangency; if they coincide then the configurations of intersection on the two sides of the stratum are still indistinguishable. We will denote by $H^1$ this type of self intersection. As explained in Remark \pr{r1}
below, this local symmetry of the $H^1$ 
configuration implies a global one-sidedness of the $H^1$ strata in the space of immersions.

If on the other hand the orientation of the two sheets disagree at time of tangency, we may distinguish the two sides of the stratum; we will consider as positive, the side for which the region in $\E$ between
the two sheets is on the preferred side of both of them. 
We denote by $H^2$ this type of self intersection.

Type $T$: This is similar to type $E$ in the sense that the two sides of the stratum are distinguishable even in the un-oriented setting, by the additional sphere which appears on one side and which we will consider as the positive side ($\lm<0$ in the formula above). This time the 2-sphere is composed of three 2-cells.
With orientation on the surface we distinguish four types of $T$ intersection, which we name $T^a$ where $0\leq a \leq 3$ is (as for $E$) the number of 2-cells
for which the little 3-cell bounded by the 2-sphere lies on their non-preferred side in $\E$.

Type $Q$: Here as for $H$, the two sides may not be distinguished without orientation; a sphere with four faces (bounding a simplex)
is formed on both sides of the stratum.
If the surface is oriented, then the numbers of faces for which the simplex is not on the preferred side in $\E$, are complementary; i.e. if it is $a$ on one side of the stratum, then it is $4-a$ on the other; so they appear in pairs $(4,0)$, $(3,1)$, $(2,2)$; 
we name the configurations, respectively $Q^4$, $Q^3$, $Q^2$. For $Q^4$ and $Q^3$ the two sides of the stratum may be distinguished and we choose the positive side to be the side where a larger number of 2-cells have the simplex on their non-preferred side.
On the other hand, for $Q^2$ the two sides are indistinguishable, and as for $H^1$, the problem is global; the $Q^2$ strata are one-sided in the space of immersions (Remark \pr{r1}).

\section{Finite Order Invariants}\label{B}

We fix once and for all a closed oriented surface $F$ and a regular homotopy class $\A$ of immersions of $F$ into $\E$.
We denote by $I_n\su \A$ ($n\geq 0$) the space of all immersions in $\A$ which have precisely $n$  CE points (the self intersection being elsewhere stable).
In particular, $I_0$ is the space of all stable immersions in $\A$.

Let $\G$ be any Abelian group and let $f:I_0\to\G$ be an invariant, i.e. a function which is constant on each connected component of $I_0$.
Given an immersion $i\in I_n$ we make an arbitrary choice of co-orientation for the configuration of intersection, at each of the $n$ CE points of $i$, and define $f^{TCO}(i)$ (Where ``TCO'' stands for ``Temporary Co-Orientation) as follows:
Let $i_1,\dots,i_{2^n}$ be the $2^n$ immersions in $I_0$ obtained from $i$ by slightly deforming it in the $2^n$ possible ways.  We define:
$$f^{TCO}(i)=\sum_{m=1}^{2^n} \epsilon^m_1\cdots\epsilon^m_n f(i_m) $$
where $\epsilon^m_k$ is $1$ or $-1$ according to whether in order to obtain
$i_m$ we deformed the configuration of $i$
at $p_k$ positively or negatively, according to the temporary 
co-orientation chosen for $i$ at $p_k$. 
The following is clear:
\begin{lemma}\label{l0}
If $i\in I_n$ and $TCO_1$, $TCO_2$ are two temporary co-orientations for $i$ which differ at precisely one CE of $i$ then $f^{TCO_1}(i) = -f^{TCO_2}(i)$.
\end{lemma}
By Lemma \pr{l0},
the statement $f^{TCO}(i)=0$ is independent of the temporary co-orientation and we may simply write $f(i)=0$.

\begin{dfn}\label{d3}
An invariant $f:I_0\to\G$ will be called \emph{of finite order} if 
there is an $n$ such that $f(i)=0$ for all $i\in I_{n+1}$.
The minimal such $n$ will be called the \emph{order} of $f$.
\end{dfn}

Let $i:F\to\E$ be an immersion having a CE located at $p\in\E$. As in [N2], we define
the degree $d_p(i)$ of the CE at $p$ as follows: Let $B$ be a tiny ball in $\E$
centered at $p$.
$i^{-1}(B)$ is a union of some (two, three or four) disks in $F$ which pass $p$.
Let $\hat{i}:F\to \pa B$ be the map obtained from $i$ as follows:
On $F-i^{-1}(B)$ we define $\hat{i}$ by radial projection (centered at $p$).
Now, if $D$ is one of the disks in
$i^{-1}(B)$ then $i(D)$ cuts $\pa B$ into two hemispheres; $\hat{i}$ is defined to
map $D$ onto the hemisphere which lies on the preferred side of $i(D)$ in $\E$. 
(Recall that $F$ is oriented.)
Finally we define $d_p(i)$ as the degree of the map $\hat{i}:F\to \pa B$ (the orientation on $\pa B$ being that induced to it from $B$.)

Let $C_p(i)$ be the expression $R^a_m$ where $R^a$
is the configuration of the CE of $i$ at $p$ (one of our twelve symbols
e.g. $E^0$) and $m=d_p(i)$. 
Let $\C_n$ denote the set of all
\emph{un-ordered} $n$-tuples of such expressions $R^a_m$.
Finally, we define a map $C:I_n \to \C_n$ as follows:
If $i\in I_n$ with CEs
located at $p_1,\dots,p_n\in\E$, then we define $C(i)\in \C_n$ to be 
the un-ordered n-tuple $[C_{p_1}(i),\dots, C_{p_n}(i)]$. 

\begin{dfn}\label{d1}
\begin{enumerate}
\item A regular homotopy $H_t:F\to\E$ will be called \emph{of type A} if it is of the form $H_t = U_t \circ i \circ V_t$ where $U_t:\E\to\E$ and
$V_t:F\to F$ are isotopies. 
\item A regular homotopy $H_t:F\to\E$ between immersions $i,j\in I_n$
will be called \emph{of type B} if it is of the following form:
If $B_1,\dots,B_n\su\E$ are little balls centered at the $n$ CE points
of $i$ and $U=i^{-1}(\bigcup_k B_k)$ then
$H_t$ fixes $U$ and moves $F-U$ within $\E - \bigcup B_k$.
\item Two immersions $i,j\in I_n$ will be called \emph{AB equivalent} if 
there is a regular homotopy $H_t$ between $i$ and $j$ which is alternatingly 
of type A and B.
Such a regular homotopy will be called an \emph{AB equivalence}.
\end{enumerate}
\end{dfn}

A proof of the following proposition appears in [N2] for immersions including
only quadruple points; the proof for the general case is identical; we include it
here for completeness.

\begin{prop}\label{p1}
Let  $i,j\in I_n$, then $i$ and $j$ are AB equivalent iff $C(i)=C(j)$.
\end{prop}

\begin{pf}
If $i$ and $j$ are AB equivalent then clearly $C(i)=C(j)$. For the converse, assume $C(i)=C(j)$. One can order the CEs of $i$ and $j$, respectively $p'_1,\dots,p'_n$ and $p_1,\dots,p_n$, such that $C_{p'_k}(i)=C_{p_k}(j)$, $k=1,\dots,n$. This means in particular, that if $B'_1,\dots,B'_n$ and $B_1,\dots,B_n$ are neighborhoods of the $p'_k$s and $p_k$s respectively, then for each $k$ there is an orientation preserving diffeomorphism from $B'_k$ to $B_k$ which takes each sheet of $i(F)\cap B'_k$ orientation 
preservingly onto the corresponding sheet of $j(F)\cap B_k$.
These diffeomorphisms may all be realized by one ambient isotopy $U_t:\E\to\E$. There is then an isotopy $V_t:F\to  F$ such that the
final immersion $i'$ of the regular homotopy $U_t\circ i \circ V_t$
satisfies that $i'$ and $j$ have the same $n$ CE points $p_1,\dots,p_n\in\E$, ${i'}^{-1}(\bigcup_k  B_k) = j^{-1}(\bigcup_k B_k)$ 
which we name $U$
and $i'|_U = j|_U$.
Also $d_{p_k}(i') = d_{p_k}(j)$ for $k=1,\dots,n$.
Now $U$ is a union of some disks $D_1,\dots,D_r$. We construct the following
handle decomposition of $F$. $D_1,\dots D_r$ will be the 0-handles. If $g$ is the
genus of $F$ we will have 1-handles  $h_1,\dots,h_{2g+r-1}$
as follows: $h_1,\dots,h_{2g}$ will
each have both ends glued to $D_1$ such that $D_1$ with $h_1,\dots,h_{2g}$ will
decompose $F$ in the standard way.
Then for $k=1,\dots r-1$, $h_{2g+k}$ will have one end glued to $D_k$ and the other to $D_{k+1}$.
The complement of the
0- and 1-handles is one disk which will be the unique 2-handle.
We will now construct a regular homotopy of the form $i'\circ V'_t$
($V'_t:F\to F$ an isotopy)
from $i'$ to an immersion $i''$ which will have the property that the restrictions of $i''$ and $j$ to all 1-handles, are regularly homotopic
keeping all 0-handles fixed.
Since $i'$ and $j$ are
regularly homotopic (recall $i,j\in I_n \su \A$),
this is already true for $h_1,\dots,h_{2g}$.
Now take $h_{2g+1}$. If $i'|_{h_{2g+1}}$ and $j|_{h_{2g+1}}$ are not regularly homotopic keeping $D_1$ and $D_2$ fixed, then $V'_t$ performs one full rotation
of $D_2$, creating a Dehn twist in a thin annulus around $D_2$ in $F$. $h_{2g+1}$ will now satisfy the needed property. Note also that this rotation of $D_2$ moves only $h_{2g+1}$ and $h_{2g+2}$, keeping all other 0- and 1-handles
fixed. We continue this way along the chain of 1-handles, rotating $D_{k+1}$ if necessary for the sake of $h_{2g+k}$.
For $k<r-1$ this will also move  $h_{2g+k+1}$,
but we never need to move 1-handles that have previously been taken care of.
Also $d_{p_k}(i'')=d_{p_k}(i)=d_{p_k}(j)$ for all $k=1,\dots, n$.
We now perform a regular homotopy $H_t$ on
the union of  0- and 1-handles which fixes the 0-handles,
and regularly homotopes each 1-handle $h$, from $i''|_{h}$ to
$j|_{h}$, avoiding $\bigcup_k B_k$. This is possible by the construction
of $i''$.
Denote our 2-handle by $D$. So far we have constructed $H_t$ only on $F-D$. By means of [S], $H_t$ may be extended to $D$, still avoiding $\bigcup_k B_k$, arriving at an immersion $i'''$.  And so, we are left with regularly homotoping $i'''|_D$ to $j|_{D}$ (relative $\pa D$). Since $d_{p_k}(i''')=d_{p_k}(j)$ for all $k=1,\dots,n$, these maps are homotopic in $\E-\bigcup_k B_k$. It then follows from the Smale-Hirsch Theorem ([H]), that they are also
\emph{regularly} homotopic in $\E-\bigcup_k B_k$ (since the obstruction
to that would lie in $\pi_2(SO_3)=0$).

The regular homotopy from $i$ to $i''$ was of type A, and that from $i''$ to $j$ was of type B.
\end{pf}

\begin{prop}\label{p2}
Let $f:I_0\to\G$ be an invariant of order $n$. Let $i\in I_n$ be an immersion 
with CEs $p_1,\dots,p_n$ and assume the configuration of $i$ at $p_1$ is of type $H^1$ or $Q^2$. Then $f^{TCO}(i) = -f^{TCO}(i)$ in $\G$.
It follows by Lemma \pr{l0} that in this case $f^{TCO}(i)$ is 
independent of the temporary co-orientation.
\end{prop}

\begin{pf}
For $k=1,\dots,n$ let $B_k$ be a small neighborhood of $p_k$ in $\E$.
Since the CE at $p_1$ is of type $H^1$ or $Q^2$,
there is an orientation preserving diffeomorphism from $B_1$ to itself which 
maps $i(F)\cap B_1$ onto itself, permuting the sheets, being orientation  preserving on each sheet, 
but reversing the co-orientation of the configuration at $p_1$.
We use this self diffeomorphism of $B_1$ and the identity map on $B_2,\dots,B_n$, in the proof of Proposition \pr{p1}, getting an AB equivalence $H_t$ from $i$ \emph{to itself} which reverses the co-orientation of the configuration in $B_1$ and fixes $F$ in $B_2,\dots,B_n$. We choose a temporary co-orientation for $i$ at $p_1,\dots,p_n$ and we carry it along $H_t$. If at some time $t_0$, $H_t$ passes through an $n+1$th CE at $p\in\E$, we choose an arbitrary co-orientation for $H_{t_0}$ at $p$. Together with the co-orientations we are carrying from $i$ this gives a temporary co-orientation for the $n+1$
CEs of $H_{t_0}$. 
Since $f$ is of order $n$, $f^{TCO}(H_{t_0})=0$ which means that 
$f^{TCO}(H_{t_0-\e})=f^{TCO}(H_{t_0+\e})$ where the TCOs at $t_0-\e$ and $t_0+\e$ are those carried from $i$ along $H_t$. Finally we arrive back at $i$ but with opposite co-orientation than the original one at $p_1$ and the same co-orientation at $p_2,\dots,p_n$. By Lemma \pr{l0} we get 
$f^{TCO}(i)=-f^{TCO}(i)$.
\end{pf}

\begin{remark}\label{r1}
We have seen in the proof of Proposition \pr{p2} that if $i\in I_n$  
has a CE of type $H^1$ or $Q^2$ then there is an AB equivalence from $i$ to itself such that if we follow that CE along the AB equivalence, then
it returns to itself but with opposite co-orientation.  
If one seeks a globally defined co-orientation, it is a minimal requirement that it be continuously defined along AB equivalences; so we see that such globally defined co-orientation does not exist for CEs of type $H^1$ and $Q^2$. 
\end{remark}

If $f$ is an invariant of order $n$ then Proposition \pr{p2} enables us to extend $f$ to $I_n$. We do it as follows: We choose once and for all a \emph{permanent} 
co-orientation for the ten configurations which allow it; 
in fact we choose those co-orientations given in Section \pr{A} above.  
Now if $i\in I_n$ and at least one of the CEs of $i$ is of configuration $H^1$ or $Q^2$ then by Proposition \pr{p2}
$f(i)$ is well defined, independent of a TCO.
If all $n$ CEs of $i$ are not of configuration $H^1$ and $Q^2$ then we define $f(i)$ using our permanent co-orientation for each of the CEs.
We will assume from now on without mention that any $f$ of order $n$ is extended to $I_n$ in this way.
(Note that if $f$ is of order $n$ then we are not extending $f$ to $I_k$ for
$0<k<n$). 

\begin{prop}\label{p3}
Let $f$ be an invariant of order $n$ and $i,j\in I_n$. 
If $C(i)=C(j)$ then $f(i)=f(j)$.
\end{prop}

\begin{pf}
By Proposition \pr{p1} $i$ and $j$ are AB equivalent.
As in the proof of Proposition \pr{p2}, $f$ is unchanged whenever we pass an $n+1$th CE  and so $f(i)=f(j)$. 
\end{pf}

By Proposition \pr{p3}
(and since $C:I_n\to\C_n$ is clearly surjective), 
any order $n$ invariant $f$ induces a well defined
function $u(f):\C_n\to\G$; if $u(f) = u(g)$ then $f$ and $g$ differ by an invariant of order at most $n-1$.
And so if $V_n$ denotes the space of all invariants on $\A$
of order at most $n$ then
$f\mapsto u(f)$ induces an injection $u:V_n / V_{n-1} \to \C_n^*$ 
where $\C_n^*$ is the space of all function from $\C_n$ to $\G$.
The purpose of the next section is to demonstrate
a space $\Delta_n \su \C_n^*$ which contains 
the image of $u$.
$\Delta_n$ will be a subspace of $\C_n^*$ which is determined by 
the restriction of Proposition \pr{p2} above
and by relations obtained by looking at local 
2-parameter families of immersions.

Comparing to knot theory,
$\C_n$ is analogous to the set of all chord diagrams of order $n$; the analogy is made clear by Propositions \pr{p1} and \pr{p3} above.
$\Delta_n$ which we define below, will be analogous to the
space of functions on chord diagrams which satisfy the 1-term and 4-term
relations.

In Section \pr{D} we will show that
$u:V_1/V_0\to \Delta_1$ is surjective; by this we classify all order one
invariants (Theorem \pr{t1}).
We will show in Section \pr{E} that for $\G={\Bbb Z}/2$ and $n>1$,
$u:V_n/V_{n-1}\to\Delta_n$ is \emph{not} surjective; this we do by demonstrating one particular function in $\Delta_n$ which is not attained by $u$.
\begin{ques}\label{qu}
What is the image of $u$ for $n>1$?
(given $F$, $\A$ and $\G$.)
\end{ques}

\section{Local analysis}\label{C}

Let $i\in \A$ be an immersion with a self intersection of local 
codim 2 at $p_1$ and 
$n-1$ additional self-intersections of local codim 1 (i.e. CEs) at $p_2,\dots,p_n$.
We look at a 2-parameter family of immersions which moves $F$ only in a neighborhood of $p_1$, such that the immersion $i$ corresponds to 
parameters $(0,0)$ and such that this 2-parameter family is transverse to the local codim 2 stratum at $i$.
In this 2-parameter family of immersions we look at a loop which circles the point of intersection with the codim 2 strata.
(This corresponds to a circle around the origin in the parameter plane.)
This circle crosses the local codim 1 strata some $r$ times.
Between each two intersections we have an immersion in $I_{n-1}$ 
with the same $n-1$ CEs, at $p_2,\dots,p_n$.
At each intersection with the local codim 1 strata, an $n$th CE is added, 
obtaining an immersion in $I_n$. Let $i_1,\dots,i_r$ be the $r$ immersions in $I_n$ so obtained and let $\e_k$, $k=1,\dots,r$ be $1$ or $-1$ according to whether we are passing the $n$th CE of $i_k$ in the direction of its 
permanent co-orientation, if it has one; if the CE is of type $H^1$ or $Q^2$ then $\e_k$ is arbitrarily chosen. 
Now let $f:I_0\to\G$ be an invariant of order $n$, 
then it is easy to see that 
$\sum_{k=1}^r \e_k f(i_k) = 0$. 
Looking now at $u:V_n / V_{n-1} \to \C_n^*$  we obtain relations that must be satisfied by a function in $\C_n^*$ in order for it to lie in the image of $u$.
In this section we will find all relations on $\C_n^*$ obtained in this way.
The relations on a $g\in\C_n^*$ will be written as relations on the symbols $R^a_m$, e.g.
$0 = T^a_m - T^{3-a}_m$ will stand for the set of all relations of the form
$0 = g([T^a_m, R^{a_2}_{d_2},\dots,R^{a_n}_{d_n}]) - 
g([T^{3-a}_m, R^{a_2}_{d_2},\dots,R^{a_n}_{d_n}])$ with arbitrary $R^{a_2}_{d_2},\dots,R^{a_n}_{d_n}$.

It will be convenient for the analysis of this section
to extend the set $\C_n$ to a set $\widetilde{\C_n}$
by also allowing the symbols $H^0$, $Q^0$, $Q^1$. 
The symbol $H^0$ will represent the same configuration as $H^2$ only with opposite co-orientation. Similarly $Q^a$ ($a=0,1$) will represent the same configuration as $Q^{4-a}$ only with opposite 
co-orientation. 
Note that our 
co-orientation for the $H$ and $Q$ configurations 
under the different labelings
may be given the following characterization: The positive side of $H^a$ is the side where $a$ of the two sheets involved have the region between the sheets on their preferred side.
The positive side of $Q^a$ is the side where $a$ of the four faces of the simplex have the simplex on their non-preferred side.
$\C_n^*$ will now be identified with the space of functions on $\widetilde{\C_n}$ 
which satisfy the relations $H^0_m=-H^2_m$, $Q^0_m=-Q^4_m$, $Q^1_m = -Q^3_m$ 
(these represent relations in the above sense).

We now need to look at local 2-parameter families of immersions which are transverse to the local codim 2 strata.
These may be divided into six types (see [HK]) which we will name by the types of CEs appearing in them:
$EH$, $TT$, $ET$, $HT$, $TQ$, $QQ$. In the notation of [HK] they are respectively: $A_0^2|A_2$, $A_0^3|A_2$, $(A_0^2|A_1^+)(A_0)$, 
$(A_0^2|A_1^-)(A_0)$, $(A_0^3|A_1)(A_0)$, $A_0^5$.

For  each of the first five types we give the following: 1. Formula for a local representative. 2. Sketch of the configuration for some value $(\lm_1,\lm_2)$ of the parameters (not $(0,0)$). 3. Diagram of the 2 dimensional parameter space, where intersection with the codim 1 strata is depicted, including their 
co-orientations (this is called a bifurcation diagram). 
4. The relation arising. 
Note that Proposition \pr{p2} takes care of the cases when there is actually no co-orientation.

For these five types,
the bifurcation diagram is obtained from the sketch and formula in a straight forward manner, by following the required loop of immersions. Whenever the plane $x=0$ appears in a
configuration below, we assume 
(by rotating the configuration if necessary) 
that its preferred side is $x>0$. 
The integer $m$ in terms of which the degrees of the CEs are given,
is the degree of the central cod 2 immersion at its cod 2 selfintersection. (We have originally defined degree only for CEs, but the same definition applies to any self-intersection.)
We then go on to type $QQ$; it requires special analysis which will be done in detail.

\begin{figure}[h]
\scalebox{0.6}{\includegraphics{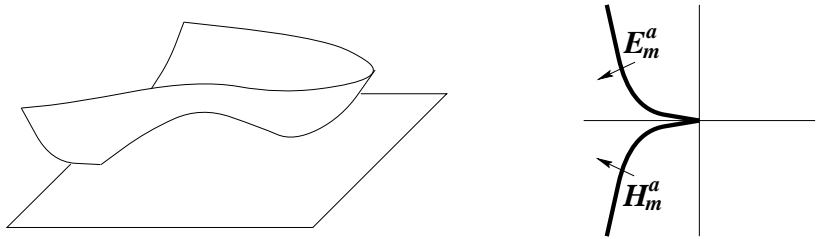}}
\caption{$EH$ configuration}\label{feh}
\end{figure}

$EH$: \ \ $z=0$, \ \ $z=y^2 + x^3+\lm_1 x + \lm_2$.
\begin{equation}\label{eeh}
0 = E^a_m - H^a_m
\end{equation}

\begin{figure}[h]
\scalebox{0.6}{\includegraphics{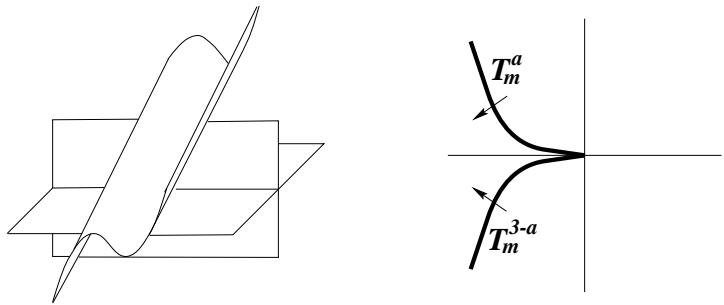}}
\caption{$TT$ configuration}\label{ftt}
\end{figure}

$TT$: \ \ $z=0$, \ \ $y=0$, \ \ $z=y+x^3+\lm_1 x  + \lm_2$.
\begin{equation}\label{ett}
0 = T^a_m - T^{3-a}_m
\end{equation}

\begin{figure}[h]
\scalebox{0.6}{\includegraphics{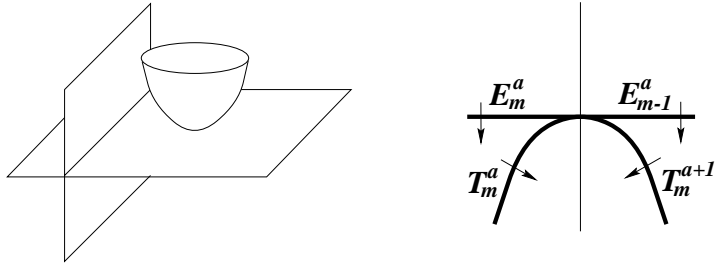}}
\caption{$ET$ configuration}\label{fet}
\end{figure}

$ET$: \ \ $z=0$, \ \ $x=0$, \ \ $z=(x-\lm_1)^2 + y^2 + \lm_2$.
\begin{equation}\label{eet}
0 = T^a_m - T^{a+1}_m - E^a_{m-1} + E^a_m
\end{equation}

\begin{figure}[h]
\scalebox{0.6}{\includegraphics{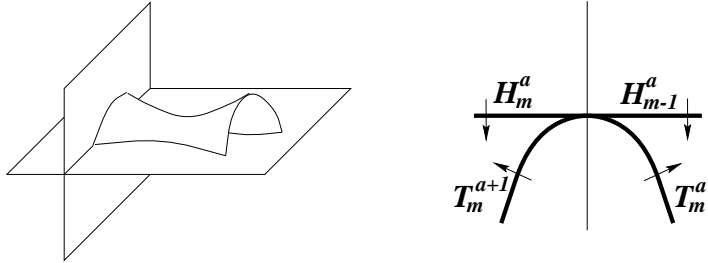}}
\caption{$HT$ configuration}\label{fht}
\end{figure}

$HT$: \ \ $z=0$, \ \ $x=0$, \ \ $z=(x-\lm_1)^2 - y^2 + \lm_2$.
\begin{equation}\label{eht}
0 = -T^{a+1}_m + T^a_m - H^a_{m-1} + H^a_m
\end{equation}

\begin{figure}[h]
\scalebox{0.6}{\includegraphics{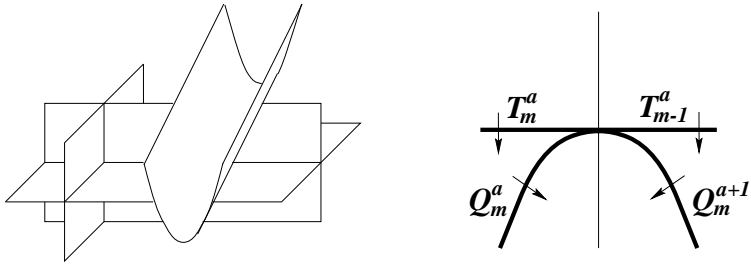}}
\caption{$TQ$ configuration}\label{ftq}
\end{figure}

$TQ$: \ \ $z=0$, \ \ $y=0$, \ \ $x=0$, \ \ $z=y+(x-\lm_1)^2 + \lm_2$
\begin{equation}\label{etq}
0 = Q^a_m - Q^{a+1}_m - T^a_{m-1} + T^a_m
\end{equation}

$QQ$: This configuration is a quintuple point, i.e. five sheets passing through a point, each three of which intersect generically.
We construct a 2-parameter family of quintuplets of oriented planes 
which represents a local 2-parameter family of immersions which is transverse to the local quintuple point stratum.
Let $P_1,\dots,P_5$ be five oriented planes through the origin in $\E$ which are in general position, i.e. no three of them intersect in a line.
For $k=1,\dd,5$ let $u_k\in\E$ be the unit vector which is perpendicular to $P_k$ and pointing into the preferred side of $P_k$ in $\E$; any three of the vectors
$u_1,\dd,u_5$ are independent.
A vector $(\lm_1,\dd,\lm_5)\in\R^5$ will represent the quintuplet of planes
$P_1^{\lm_1},\dots,P_5^{\lm_5}$ in $\E$ where 
$P_k^{\lm_k}=\{ x\in\E : x\cdot u_k = \lm_k \}$.
(In particular, for each $k$: $P^0_k = P_k$).

If $\lm\in\R$ and $v\in\E$ then $P_k^{\lm+v\cdot u_k}$ is the translate by the vector $v$ of the plane $P_k^\lm$. Let $V\su\R^5$ be the
3-dimensional subspace defined by
$V=\{ (v\cdot u_1, \dd , v\cdot u_5) : v\in \E \}$ then 
$(\lm_1,\dots,\lm_5)\in V$ iff $P_1^{\lm_1},\dots,P_5^{\lm_5}$ all meet at a point.
If $U$ is a direct summand of $V$ in $\R^5$  then $U$ represents all configurations of quintuplets 
of planes which are respectively parallel to $P_1,\dd,P_5$, up to a common translation of the five planes. Such a $U$ is then  a representative local 2-parameter family of immersions which is transverse to the local codim 2 stratum of quintuple points.
We make the choice $U=V^\perp$. If $A$ is the $5\times 3$ matrix whose rows are  $u_1 , \dd , u_5$
then the columns of $A$ span $V$ and so
$U=V^\perp$ is the left kernel of $A$
i.e. $U=\{ (\lm_1,\dd,\lm_5)  : \lm_1 u_1 +\cdots + \lm_5 u_5 = 0 \}$.

We now find the points in $U$ which represent configurations involving a 
CE, which in this case must be a quadruple point. Let $l_k = U\cap ( V + \R \e_k)$ where $\e_1,\dd,\e_5$ is the standard basis of $\R^5$ (i.e. $\e_1=(1,0,0,0,0)$ etc.);
then $l_k$ represents all configurations in $U$ where the planes
$\{P_j : j\neq k \}$ meet at a quadruple point; indeed, an element of $l_k$ is
an element of $U$ which is
obtained from an element of $V$ i.e. a quintuple point,
by adding some $r\e_k$ i.e. pushing away the plane $P_k$.

For a given $l_k$ we would like to determine the configuration of the quadruple point represented by it. Take say $l_5$. We must look at the configuration of the four planes $P_1^{\lm_1},\dd,P_4^{\lm_4}$ given by points 
$(\lm_1,\dots,\lm_5)$ on the two sides of $l_5$ in $U$,
and see how many of the planes $P_1^{\lm_1},\dd,P_4^{\lm_4}$ 
have the simplex created by them, on their non-preferred side.
It is enough to check one point on each side of $l_5$ in $U$ and as we shall see, it will be most convenient to look at the points of
$l_5^\perp$ (here $\perp$ means the orthogonal complement in $U$).
Now $l_5^\perp =  \e_5^\perp \cap U$ (the $\perp$ on the right is
the orthogonal complement in $\R^5$);
i.e. $l_5^\perp
=  \{ (\lm_1,\dd ,\lm_4, 0 ) \in \R^5 : \lm_1 u_1 + \cdots + \lm_4 u_4 = 0 \}$.
For such points we can determine the configuration of the simplex
using the following lemma:

\begin{lemma}\label{l1}
Let $v_1,\dd,v_4\in\E$ be four vectors such that any three of them are independent and $v_1 + \cdots + v_4 =0$.
If $\mu_1,\dd,\mu_4$ are \emph{positive} real numbers
then the origin $0\in \E$ lies in the interior of the simplex determined by the four planes
$\{x\in\E : x\cdot v_k = \mu_k \}$.
\end{lemma}

\begin{pf}
Since $0\cdot v_k = 0 < \mu_k$, the domain determined by the four planes
in which the origin lies is:
$D=\{ x\in\E : x\cdot v_k \leq \mu_k \ \ \hbox{for all}  \ \ 1\leq k \leq 4 \}$.
This is a convex domain in $\E$. If it is not the simplex determined by the four planes then it is unbounded and so there
is a ray based at 0 which is contained in $D$, i.e. there is a vector $v\neq 0$ such that $rv \cdot v_k \leq \mu_k$ for every $1\leq k \leq 4$ and any $r>0$.
It follows that for each $k$: $v \cdot v_k \leq 0$.
If for some $k$, $v\cdot v_k < 0$ then $v \cdot (v_1 + \cdots + v_4) < 0$
contradicting $v_1 + \cdots + v_4 =0$. So $v\cdot v_k =0$ for all $k$, contradicting the fact that any three of $v_1,\dd,v_4$ are independent.
\end{pf}

Back to our $l_5^\perp$,
we use Lemma \pr{l1} with $v_k = \lm_k u_k$ and $\mu_k = (\lm_k)^2$, $k=1,\dd,4$;
obtaining that 0 lies in the interior of the simplex determined by the equations
$x\cdot \lm_k u_k = (\lm_k)^2$ ($k=1,\dd,4$) which is the same as
$x\cdot u_k = \lm_k$ i.e. the planes $P_1^{\lm_1}, \dd , P_4^{\lm_4}$.

Now, the vector $u_k$, when based at $P_k^{\lm_k}$, points into the preferred
side of $P_k^{\lm_k}$ in $\E$. On the other hand, 
$u_k$ points away from the side where the origin lies, which is the side 
where the simplex lies, iff $\lm_k > 0$.
We conclude that if
$(\lm_1,\dd,\lm_4,0)\in l_5^\perp$ then the number of faces of the simplex determined by $P_1^{\lm_1},\dd,P_4^{\lm_4}$ 
which have the simplex on their non-preferred side, is precisely the number of positive numbers among $\lm_1,\dd,\lm_4$. (Note that $\lm_1,\dd,\lm_4$ are all non-zero since each three $u_k$s are independent).

The origin of $U$ splits $l_5^\perp$ into two half lines; clearly the number $p$ of plus signs is constant on such a half-line and is $4-p$ on the other half line.
So for given quintuple point, our task is to find the number $p$ of plus signs
in each of the ten half lines of the $l_k^\perp$s.
Now $l_k^\perp$ has $\lm_k=0$, so it partitions $U$ into the domains where
$\lm_k>0$ and $\lm_k<0$; we thus use it to determine the sign of $\lm_k$ in the half-lines of $l_j^\perp$ for $j\neq k$. 
Examples of such analysis appear in Figure \pr{fqq}, which we now explain:
Each of the short thick lines in a diagram represents an $l_k^\perp$. The arrows on each $l_k^\perp$ point to the 
side of it in $U$ where $\lm_k$ is positive. 
To determine the number $p$ of positive $\lm_j$s corresponding to a given half line of $l_k^\perp$ we need to count for how many $l_j^\perp$s  ($j\neq k$)
this half line lies on their positive side, i.e. the side designated by the arrows.
This is the number appearing in the circle located on the given half $l_k^\perp$ in the diagram (at the tip of the short thick line).
Finally the longer thinner lines are the $l_k$s themselves (each drawn perpendicular to the corresponding $l_k^\perp$).
The pair of numbers at each tip of $l_k$ is simply copied
from the corresponding sides of $l_k^\perp$. When passing an $l_k$, this pair of numbers (which is of the form $p, 4-p$) tells us the type of quadruple point we are passing, and the co-orientation with which we are passing it, namely, 
the side chosen by our permanent co-orientation is
the side where the larger number of the pair appears.
(As before, Proposition \pr{p2} takes care of the case when there is no
co-orientation.)

\begin{figure}[ht]
\scalebox{0.7}{\includegraphics{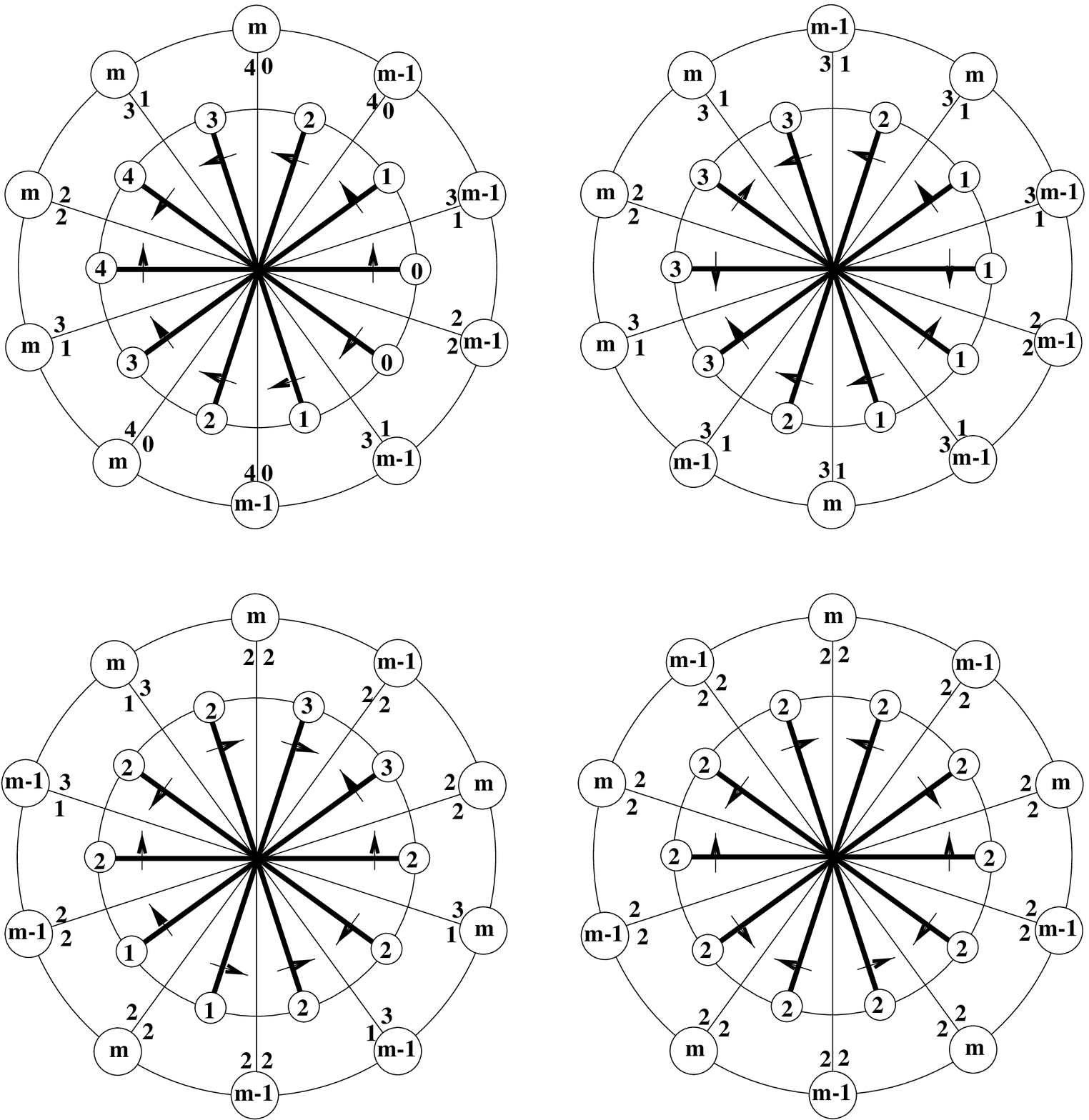}}
\caption{$QQ$ configuration}\label{fqq}
\end{figure}

Finally we need to determine the degrees $d_p(i)$ at each of the ten quadruple points. Let $m$ be the degree of the quintuple point and look say at $l_1$.
We claim that the half of $l_1$ which has $\lm_1>0$ represents quadruple points with degree $m$ whereas the half with $\lm_1<0$ represents quadruple points with degree $m-1$. 
To establish this we need to show that a quadruple point represented by a point in the half of $l_1$ with $\lm_1>0$ is obtained from a translate of our quintuple point by pushing $P_1$ into its preferred side i.e. the side pointed at by $u_1$. Recall that a point $(\lm_1,\dots,\lm_5)\in l_1$ is of the 
form $r\e_1 + v$ where $v$ is perpendicular to $U$ and 
represents a common translation of all five planes;
so $P_1^{\lm_1}$ is obtained by pushing $P_1$ away from the quintuple
point represented by $v$, and this push is into the preferred side of $P_1$
iff $r>0$.
It remains to notice that the sign of $\lm_1$ is the same as that of $r$,
since $(\lm_1,\dots,\lm_5)$ is the orthogonal projection of $(r,0,0,0,0)$ to $U$. The integers $m$ or $m-1$ appearing at the tips of each $l_k$ in the diagrams are the degrees.

Once we have all this information registered at the tips of the $l_k$s, we can read off the diagram the relation determined by the given quintuple point. 
We claim that for any configuration of lines and arrows as above, 
eight of the terms cancel, always leaving the same relation:
\begin{equation}\label{eqq}
Q^2_m = Q^2_{m-1}
\end{equation}
We offer two ways of seeing this: The first way is to verify that the four diagrams appearing in Figures \pr{fqq} are in fact all essentially distinct ways for choosing such distribution of arrows; then explicitly write out the relation obtained by each of them. We only remark that there is no loss of generality by the fact that the diagrams are sketched with equal angles between the $l_k$s, since all that is relevant to us is the cyclic ordering of the $l_k^\perp$s and not the common cyclic ordering of all $l_1,\dd,l_5,l_1^\perp,\dd,l_5^\perp$.

A second way to see that the same relation $Q^2_m=Q^2_{m-1}$ is always obtained is as follows: Since we have five arrows pointing clockwise and five counter-clockwise, there must be two consecutive arrows which are pointing at each other.
These half $l_k^\perp$s must have the same number at their tip, since they are on the same side of any other $l_j^\perp$ and are both on the positive side of each other (i.e. the side designated by the arrow). We now look at the corresponding positive half $l_k$s. They both have degree $m$ and will both have the same number copied from $l_k^\perp$ next to them, but written on the opposite side. So their contribution to the relation is precisely the negative of each other;
the same will be true for the negative halves of these $l_k$s (only now with degrees $m-1$.)
We now erase these two lines from the diagram and argue by (a two step) induction that the formal property analogous to the one we are proving, holds for the remaining three line diagram (checking that it holds for a one line diagram.) It remains to observe that once we re-insert the two erased lines, they jointly add precisely 1 to the number at the end of all other short thick lines and so the property holds for the five line diagram.

We denote by $\Delta_n = \Delta_n(\G)$ 
the subspace of $\C_n^*$ satisfying relations \pr{eeh}-\pr{eqq} 
above and the restriction of Proposition \pr{p2}, which we may write as 
$0= 2H^1_m = 2Q^2_m$. Recall also that by definition of $\C_n^*$ when using 
the set $\widetilde{\C_n}$, we also have the relations $H^0_m=-H^2_m$ 
and $Q^a_m=-Q^{4-a}_m$ ($a=0,1$).

Now that we have obtained our relations, 
we return to our original $\C_n$ with only twelve symbols
i.e. we dispose of the redundant symbols $H^0$, $Q^0$ and $Q^1$. Let $\B\su\G$ be defined by
$\B=\{ x\in \G : 2x=0\}$. 
After some simplification, the relations defining $\Delta_n$
may finally be presented as follows:

\begin{itemize}
\item $E^2_m = - E^0_m = H^2_m$, \ \  $E^1_m = H^1_m$.
\item $T^0_m = T^3_m$, \ \  $T^1_m = T^2_m$.
\item $H^1_m = H^1_{m-1} \in \B$.
\item $Q^2_m = Q^2_{m-1} \in \B$.
\item $H^2_m - H^2_{m-1} = T^3_m - T^2_m$
\item $Q^4_m - Q^3_m = T^3_m - T^3_{m-1}$, \ \  $Q^3_m - Q^2_m = T^2_m - T^2_{m-1}$ 
\end{itemize}

\section{Order One Invariants}\label{D}

In this section we will show that the injection $u:V_1 / V_0 \to \Delta_1$ is surjective. 
Let us first give an explicit presentation of $\Delta_1$. We see from the presentation of the relations
appearing in the end of the previous section, that a function $g\in\Delta_1$ 
may be assigned arbitrary values in $\G$ for the symbols
$\{T^2_m\}_{m\in{\Bbb Z}}$, $\{T^3_m\}_{m\in{\Bbb Z}}$, $H^2_0$ 
and arbitrary values in $\B$ for the two symbols $H^1_0$, $Q^2_0$. Once this is done then the value of $g$ on all other symbols is uniquely determined; namely: 
\begin{enumerate}
\item $H^1_m = H^1_0$ for all $m$.
\item $H^2_m = H^2_0 + \sum_{k=1}^m (T^3_k - T^2_k)$ for $m\geq 0$.
\item $H^2_m = H^2_0 - \sum_{k=m+1}^0 (T^3_k - T^2_k)$  for $m<0$.
\item $E^0_m = -H^2_m$, \ \ $E^1_m = H^1_m$, \ \  $E^2_m = H^2_m$ for all $m$.
\item $T^0_m = T^3_m$, \ \ $T^1_m = T^2_m$ for all $m$.
\item $Q^2_m = Q^2_0$ for all $m$.
\item $Q^3_m = Q^2_m + T^2_m - T^2_{m-1}$ for all $m$.
\item $Q^4_m = Q^3_m + T^3_m - T^3_{m-1}$ for all $m$.
\end{enumerate}

Let $X$ denote the set of symbols
$\{T^2_m\}_{m\in{\Bbb Z}}\cup\{T^3_m\}_{m\in{\Bbb Z}}\cup\{H^2_0\}$ 
then we have obtained that 
$\Delta_1(\G)\cong\G^X \oplus \B \oplus \B$ where 
$\G^X$ denotes the group of all functions from $X$ to $\G$.
We can also present $\Delta_1(\G)$ through a universal object as follows:
We define a universal Abelian group $\G_U$ by the Abelian group presentation
$\G_U = \left< \{t^a_m\}_{a=2,3, m\in{\Bbb Z}}, h^2_0, h^1_0, q^2_0 \ | \ 
2h^1_0 = 2q^2_0 = 0 \right>$. 
Then we define the universal element $g_U\in\Delta_1(\G_U)$ by $g_U(T^a_m) = t^a_m,
g_U(H^2_0)=h^2_0, g_U(H^1_0)=h^1_0, g_U(Q^2_0)=q^2_0$ 
and the value of $g_U$ on the other symbols of $\C_1$ is determined by
formulae 1-8 above, so indeed $g_U\in\Delta_1(\G_U)$. 
Then for arbitrary Abelian group $\G$ we have 
$\Delta_1(\G) \cong Hom(\G_U , \G)$ where the isomorphism 
maps a homomorphism $\phi:\G_U \to \G$ to the function $\phi \circ g_U \in \Delta_1(\G)$.
We will show that there is an order 1 
invariant $f_U:I_0\to\G_U$ such that its extension 
to $I_1$ induces $g_U$ on $\C_1$, i.e. $u(f_U)=g_U$. It will follow that for any group $\G$, $u:V_1/V_0\to \Delta_1(\G)$ is surjective, since if $g\in\Delta_1(\G)$ and $g=\phi \circ g_U$ where $\phi\in Hom(\G_U , \G)$ 
then $u(\phi\circ f_U) = g$.

We choose a base immersion $i_0\in I_0$ once and for all. Given any $i\in I_0$ we take a generic regular homotopy from $i_0$ to $i$ i.e. a path in $\A$ from $i_0$ to $i$ transverse with respect to the global stratification. The value of $f_U$ on $i$ will be the sum with signs, of the values of $g_U$ 
on the CEs that we pass along our path from $i_0$ to $i$, where the signs are determined by whether we are passing the given CE
in the direction of its co-orientation.  
(Note that whenever there is no co-orientation then the element we are adding
is of order 2 in $\G_U$.)
We must show that this sum is independent of our choice of path, or equivalently that it is 0 along any closed path.
First we observe that it is 0 on null-homotopic paths. Indeed, by slightly deforming the null-homotopy, we may assume it too is transverse with respect to the stratification. If we then break the null-homotopy into small pieces and look at the value going around each little piece, 
then each is 0 since $g_U\in\Delta_1(\G_U)$. Note that the 
symbols encountered when going around a global codim 2 stratum which corresponds to two local codim 1 intersections occurring at distinct places,  always adds up to 0.

Once we have established that this sum in $\G_U$ is 0 for any null-homotopic loop, we have a well defined homomorphism $\pi_1(\A)\to\G_U$ and we must show that this homomorphism is 0.
We first show this for the case $F=S^2$ in which case there is only one regular homotopy class, which we still name $\A$. Now $\pi_1(\A)=\Z \oplus \Z/2$ where the unique order 2 element is represented by one full rotation of $S^2$ 
in $\E$.
Such rotation passes through no CEs at all, and so gives 0 as needed.
Let $K\su \pi_1(\A)$ be the subgroup generated by this order 2 element in $\pi_1(\A)$ and let $P=P(S^2)=\pi_1(\A) / K$; then we are left with showing that the map induced on $P$ ($\cong \Z$) is 0.

If $H_t:S^2\to\E$ is a regular homotopy, then attaching to $H_t$ one or more 
of the superscripts $\E$, $S^2$ and $T$ will denote respectively, composition from the left with the antipodal map of $\E$, composition from the right with the antipodal map of $S^2$ and reversal of time. Note that 
these three operations on $H_t$ commute with each other. 
Let $p:S^2\to \R P^2$ be the double covering and let $i:\R P^2\to\E$ be some immersion. Let $s=i\circ p$ then $s:S^2\to\E$ is an immersion satisfying 
$s(-x)=s(x)$. Let $H_t$ be a regular homotopy from the inclusion $e:S^2\to \E$ to $s$. Let $G_t = H_t *(H_t^{S^2,T})$ where $*$ denotes concatenation from left to right, then $G_t$ is a regular homotopy from $e$ to $-e$ i.e. an eversion of the sphere. 
Also $G_t^{S^2,T} = G_t$.
We now define $F_t = G_t * G_t^\E = G_t * (G_t^{S^2,T})^\E = G_t * (G_t^{S^2,\E})^T$. 
From the presentation of $F_t$ as $G_t * G_t^\E$ we can see that 
there are isotopies $U_t:\E\to\E$ and $V_t:S^2\to S^2$ such that 
$U_t \circ F_t \circ V_t$ satisfies the conditions of [N1] Proposition 2.1.
It follows that $F_t$ represents an odd power of the generator of $P$.
We also have $F_t=G_t * (G_t^{S^2,\E})^T$ and let $G'_t$ be a generic 
regular homotopy which is a slight perturbation of $G_t$ 
($G_t$ was highly non-generic), then 
$F'_t=G'_t * ({G'}_t^{S^2,\E})^T$ will still represent an odd power of the generator of $P$ (though it will not be precisely equal to $G'_t * {G'}_t^\E$). 
We now observe that reversing the orientations of both $S^2$ and $\E$ simultaneously, preserves all our definitions of types, degrees and co-orientations of CEs and so ${G'}_t^{S^2,\E}$ produces precisely the same element in $\G_U$ as $G'_t$ and so $({G'}_t^{S^2,\E})^T$ produces the negative of that element. So we get that the value on the loop $F'_t$ is 0.
Since $F'_t$ represents an odd power of the generator of $P$ and since in $\G_U$ 
an odd multiple of any non-zero element is still non-zero, we get that $g_U$ is 0 on the generator of $P$ and so on all $P$. 
This completed the proof for the case $F=S^2$.

Now let $F$ be of higher genus and let $\A$ be a given regular homotopy class.
Again $\pi_1(\A) \cong \Z \oplus \Z/2$. Again our claim is clear for 
$K\su \pi_1(\A)$ which is defined as above and so again it is enough to look at the map induced on $P=P(F,\A)=\pi_1(\A) / K$.

Looking back for a moment at $S^2$ let $H_t$ be a loop representing a generator of $P(S^2)$ and we may assume $H_t$ fixes a given disc $D\su S^2$. Choose a point $x$ on a ray perpendicular to the fixed image of $D$ in $\E$, $x$ being far enough so that the image of $H_t$ never passes $x$. Now change the constant embedding of $D$ to be one with a thin ``thorn'' pulled out of $D$ and embedded along this ray, reaching $x$.
The new loop $H'_t$ obtained in this way
also represents the generator of $P(S^2)$
and  has the property that a given point in the fixed disc $D$ (namely, the tip of the thorn) does not participate in any self intersections throughout $H'_t$. 
Since we have already proved our result for $S^2$, we know that the symbols encountered along the loop $H'_t$ add up to 0 in the group $G_U$.
Now take a tiny immersion of $F$ in our regular homotopy class $\A$,
located near the tip of our fixed thorn, and connect sum it with $H'_t$ for all $t$, obtaining a loop $F_t:F\to\E$ in $\A$. 
It follows from the proof of [N1] Theorem 3.4
that this loop $F_t:F\to\E$ represents a generator of $P(F,\A)$.
Since the tip of the thorn was far away from any self intersections occurring
in $H'_t$, the tiny immersion will not participate in any CEs occurring in $F_t$ and so will not add to the sum of symbols attained by $H'_t$, which we know is 0 in $G_U$. 
We conclude that $f_U$ is well defined; it is evident that $f_U$ is of order 1 and
$u(f_U)=g_U$.
This completes the proof of our main result:

\begin{thm}\label{t1}
For any closed orientable surface $F$, regular homotopy class $\A$ of immersions of $F$ into $\E$ and Abelian group $\G$, the injection
$u:V_1 / V_0 \to \Delta_1$ is surjective. 
\end{thm}

\section{Higher order invariants}\label{E}

In this section we show that for $\G=\Z/2$ and $n>1$ the map
$u:V_n / V_{n-1} \to \Delta_n$ is not surjective.
We define a function $g:\C_n\to\Z/2$ as follows: 
If $x\in \C_n$ includes at least one symbol which is not of type $Q$, then $g(x)=0$ and if all symbols in $x$ are of type $Q$ then $g(x)=1$. One verifies directly that $g$ satisfies the relations appearing in the end of Section \pr{C}
and so $g\in \Delta_n(\Z_2)$. We will show that there is no order $n$ invariant $f$ such that $u(f)=g$. Let $f:I_0\to \Z/2$ be an invariant of order $n$ and assume $u(f)=g$.
Since $\G=\Z/2$, there is no need for co-orientations for the extension of $f$ and so $f$ extends to $I_k$ for any $k>0$, in particular to $I_{n-1}$.
Let $i,j\in I_{n-1}$ be two AB equivalent immersions such that all $n-1$ CEs of $i$ and $j$ are of type $Q$.
If $H_t$ is a generic AB equivalence between $i$ and $j$ then there are some finitely many times along $H_t$ for which an $n$th CE appears. By definition of $g$, the number mod 2 of such $n$th CEs which are of type $Q$ is 
$f(i)-f(j)\in\Z/2$. 
In particular if $i=j$ then the number of such occurrences of an $n$th quadruple point is 0 mod 2. 

Now take some $i\in I_{n-1}$ with $n-1$ CEs located at $p_1,\dots,p_{n-1}$ all of type $Q$ and furthermore the CE at $p_1$ is of type $Q^2$.
Let $H_t$ be an AB equivalence from $i$ to itself which fixes $F$ in neighborhoods of
$p_2,\dots,p_{n-1}$ and reverses 
the co-orientation of  the CE at $p_1$
(as in the proof of Proposition \pr{p2}). 
By the above discussion, the number mod 2 of quadruple point occurrences 
along $H_t$ is $0$.
Let $i'\in I_0$ be an immersion which is obtained from $i$ by slightly deforming $i$ in a neighborhood of each $p_k$, $k=1,\dots,n-1$.
$H_t$ induces a regular homotopy $H'_t$ 
from $i'$ to an immersion which differs from $i'$ only in a neighborhood of $p_1$ and there are 0 mod 2 quadruple point occurrences along $H'_t$, precisely those occurring during $H_t$. 
We complete $H'_t$ to a closed loop ending at $i'$ causing one more quadruple point to occur, namely that at $p_1$.
And so we have constructed a generic closed loop in $\A$ with 1 mod 2 quadruple point
occurrences. This would imply that the analogous order 1 invariant does not exist; but we have shown (Theorem \pr{t1})
that for order 1, all invariants defined in this way do exist; a contradiction. 
We conclude that 
there is no order $n$ invariant such that $u(f)=g$ and so for $\G=\Z/2$
and $n>1$ the map $u:V_n / V_{n-1} \to \Delta_n(\Z / 2)$ is not surjective.

We remark that the order 1 invariant whose existence has been used in the above proof, namely that order 1 invariant which counts the number mod 2 of quadruple points occurring along generic regular homotopies, has been studies in [N1] and [N3]. In [N1] its existence has been established for all closed surfaces (not only orientable) and also 
the existence of the analogous invariant in general 3-manifolds, 
under certain conditions. In [N3] an explicit formula has been given 
for the value of this invariant on embeddings.

We conclude with another remark. In [N2] the notion of a $q$-invariant has been defined. This is an invariant $f:I_0\to\Z/2$ such that its extension to $I_1$ satisfies that $f(i)=0$ if the CE of $i$ is not a quadruple point. 
It follows that $f(i)=0$ for any 
$i\in I_k$ ($k\geq 1$) 
which includes at least one CE which is not a quadruple point. 
If $f$ is of order $n$ then this means that $u(f)(x)=0$ for any $x\in\C_n$ which includes at least one symbol which is not of type $Q$.
It has been shown in [N2] and also follows from the relations defining $\Delta_n$ 
(end of Section \pr{C}), that for such $f$ of order $n$, $f(x)=1$ for any $x\in\C_n$ 
for which all symbols are of type $Q$. 
So what we have shown in this section is that $q$-invariants of order $n>1$ do not exist;
by this completing the work of [N2]. 
(The unique $q$-invariant of order 1 is the invariant discussed in the 
previous paragraph.)


\end{document}